\documentclass[12pt]{report}
\usepackage{amsmath}
\usepackage{amsfonts}
\usepackage{amssymb}
\usepackage{mathrsfs}
\usepackage{graphicx}
\usepackage{cancel}
\usepackage{hyperref}
\usepackage[english]{babel}
\addto{\captionsenglish}{%
  
}
\DeclareMathAlphabet{\mathpzc}{OT1}{pzc}{m}{it}

\begin{document}

\begin{center}
\textbf{\large A Generalization of Schur's Theorem}

\end{center}
\begin{center}
\textbf{Jon Henry Sanders }
\end{center}
$$$$
\begin{center}
JHS Consulting          jon\_sanders@partech.com\\

\end{center}

\chapter{}
\section{\small A GENERALIZATION OF SCHUR'S THEOREM}
\noindent
\par In 1916 Schur [1] proved the following theorem:\\
\underline{Theorem 1.}~~~~
Given a number\footnote{The word ``number'' will mean ``positive integer'' throughout this section. } there exists a number $J(t)$ such that if the set $\{1,\ldots,J(t)\}$ is divided into $t$ classes $A_1,\ldots,A_t$ then there exists numbers $a,b$  and a number $\ell, 1\leq \ell \leq t,$ such that $a,b, ~~ a+b\in A_\ell$. We prove the following generalization:\\
\underline{Theorem 2.}~~~~
Let $	P(t,n)$ be the statement, ``Given numbers $t,n$ there exists a number $J(t,n)$ such that if the set $S=\{1,\ldots,J(t,n)\}$ is divided into $t$ classes $A_1,\ldots,A_t$ there exists a sequence of (not necessarily distinct) numbers  $a_1,\ldots,a_n$ and a number $\ell,~~ 1\leq \ell\leq t$, such that $  \sum_{i\in I}a_1\in A_\ell$ for all $I\subseteq\{1,\ldots,n\}$." Then $P(t,n)$ is true for all numbers $t,n.$
\par We first prove the following corollary:\\
\underline{Corollary 1.1.} ~~~The ${a_i}$ may be assumed distinct in the above theorem.
\underline{Proof.}~~~~ Let $J^\prime(t,n)=J(t,n^4).$ We show if the set $S=\{1,\ldots,J^\prime(t,n)\}$ is divided into $t$ classes there are $n$ distinct numbers in $S$ with the desired property. For let $a_1,\ldots,a_{n^4}$ be a sequence of $n^4$ not necessarily  distinct numbers of $S$ with the desired property. Then either there are at least $n^2$ distinct numbers among them or there is a number $k$ such that $a_1=k$ for at least $n^2$ distinct numbers among them or there is a number $k$ such that $a_i=k$ for at least $n^2$ distinct values of $i$. In the first case we are through. In the second define $b_i=ik,~~ i=1,\ldots,n.$ Then $\displaystyle \sum_{i\in I}b_i, I\subseteq\{1,\ldots,n\}$ is equal to 
$\ell k,~ 1\leq\ell\leq\frac{n(n-1)}{2}\leq n^2=\displaystyle \sum_{i\in I}a_i,I^\prime\subseteq\{1,\ldots,n^4\}$ where $I^\prime$ is a set of $\ell$ (distinct) indices which satisfy $i\in I^\prime$ implies $a_i=k.$ Thus the $b_1$ have the desired property since the $a_1$ do and furthermore they are distinct.$\Box$
\par We need to make one remark before proceeding. Let $P(t,n,k)$ be the statement derived from $P(t,n)$ by substituting the set $S^\prime=\{k,2k,\ldots,J(t,n)k\}$ for the set $S=\{1,2,\ldots,J(t,n)\}$.\\
\underline{Remark.} For any given number $k,P(t,n)$ implies\footnote{Actually, $P(t,n)$ and $P(t,n,k)$ are equivalent. }  $P(t,n,k)$.\\
\underline{Proof.}~~~~ Let $A_1,\ldots,A_t$ be a division of $S^\prime$ into $t$ sets. 
Divide $S$ into $t$ sets by the rule for all $a\in S,a\in A^\prime_j$ iff $ka\in A_j$. If $P(t,n)$ is true there exists  a sequence $b_1,\ldots,b_n$ such that $\displaystyle \sum_{i\in I}b_1\in A^\prime_\ell$ for some fixed $\ell,~~1\leq\ell\leq t,$ for all $I\subseteq\{1,\ldots,n\}.$ But then the sequence $a_i=b_ik,~ i=1,\ldots,n$ and the number $\ell,~1\leq\ell\leq t,$ have the desired property for $S^\prime$ so $P(t,n,k)$ is true.$\Box$
\par We now prove an ``Iterated Ramsey Theorem'' needed for the proof of Theorem 2. Let $R(k,r,t)$ be Ramsey's function, i.e., let $R(k,r,t)$ be the smallest number such that when the $r-$subsets of any set $S$ of order $R(k,r,t)$ are divided into $t$ classes $A_1,\ldots,A_t$ there always exists a $k$-subset  $K$ of $S$ and a number $\ell,~ 1\leq\ell\leq t$ such that all the $r$-subsets of $K$ are contained in $A_\ell$. We prove\\
\underline{Lemma 1.}~~~~
Given numbers $k,r,t, ~~~~ k\geq r$ there exists a number $N(k,r,t)$ such  that if \underline{all} the non-empty subsets of a set $S$ of order $n\geq N(k,r,t)$ are divided into $t$ classes $A_1,\ldots,A_t$ then there exists a $k-$subset $K\subseteq S$ and a sequence $A_{i_1},\ldots,A_{i_r},~~1\leq i_j\leq t, ~~~~j=1,\ldots,r$, such that all the $j$-subsets of $K$ are contained in $A_{i_j}$ for $j=1,\ldots,r.$\\ 
\underline{Proof.}~~ The proof is by induction on $r$. If $r=1$ we take $N(k,l,t)>kt$. Suppose the theorem is true for $r=r_o{-1}$ and arbitrary $k,t$. We wish to prove it true for $r=r_o$ and arbitrary values of $k,t$, say $k_o,t_o$. Take $N(k_o,r_o,t_o)=N(R(k_o, r_o{-1},~t_o),~r_o,~ t_o)$. If the order of $S$ is larger than $N(k_o,r_o,t_o)$ then by the induction  hypothesis there exists a subset $K^\prime$ of $S$ of order $R(k_o,r_o,t_o)$ and a sequence $i_1,\ldots,i_{r_o{-1}}$ such that each $j$-subset of $K^\prime$ is contained in $A_{i_j},j=1,\ldots,r_o{-1}.$ But by Ramsey's Theorem there exists a $k_o-$subset $K^{\prime\prime}$ of $K^\prime$ such that all the $r_o-$subsets of $K^\prime$ are in some $A^\prime$. Then $K^{\prime\prime}$ is the desired subset of $S$ and $A_{i_1},\ldots,A_{i_{r_o{-1}}}, A^\prime$ is the required sequence of classes.$\Box$
\par We are now ready to prove Theorem 2.\\
\underline{Proof}.~~The proof is by induction on $n$. Let $P^\prime(n)$ be the statement $``P(t,n)$ is true for all possible values of t.'' Then $P^\prime(2)$ is Theorem 1. Assume $P^\prime(n_o{-1})$ is true. We prove $P^\prime(n_o)$ is true. Take
$$J(n_o,t)=N(4J(n_o{-1},t),~ 4J(n_o{-1},~t),~t).$$
Let the set $S=\{1,\ldots,J(n_o,t)\}$ be divided into $t$ classes $A_1,\ldots,A_t$. If $B=\{a_1,a_2,\ldots,a_k\}$ is any subset of $S,~ a_1<a_2<\cdots<a_k $ define $f(B)$ by $f(B)=\displaystyle \sum_{j=1}^k(-1)^{j}a_j$ if $k$ is even and  by $f(B)=\displaystyle \sum_{j=1}^k(-1)^{j+1}a_j$  if $k$ is odd.\\
Divide all the non-empty subsets of $S$ into $t$ classes $A_1^\prime,\ldots,A_t^\prime$ by the rule $B\in A_i^\prime$
 iff $f(B)\in A_i$. (Note that this is well defined since $f(B)\in S$ if $\theta\neq B\subseteq S$). By Lemma $1$\footnote{with $k=r=N.$}, there exists a subset $K$ of $S$ of order $4J(n_o^{-1},t)$ := $N$ and a sequence of classes $A_{j_1}^\prime,\ldots,A_{j_N}^\prime$ such that all the i-subsets of $K$ are contained in $A_{j_i}^\prime,i=1,\cdots,N.$ Now divide the set $S^\prime=\{4,8,\ldots,4J(n_o{-1},t)\}$ into $t$ classes $A^{\prime\prime}_1,...,A_t^{\prime\prime}$ by the rule for all $a\in S^\prime,~ a\in A_j^{\prime\prime}$ iff all the $a-$subsets of $K$ are contained in $A_j^\prime$. Then by the induction hypothesis and the preceding remark we can find a sequence $a_1,\ldots,a_{n_o{-1}}$, say $a_1\leq a_2\leq\cdots\leq a_{n_o{-1}}$, and a number $\ell,~~ 1\leq\ell\leq t$, such that $\displaystyle \sum_{i\in I}a_i\in A_\ell^{\prime\prime}$ for all $I\subseteq\{1,\ldots,n_o{-1}\}$. Since $\displaystyle \sum_{i=I}a_i\leq$ order of $K$\footnote{Since $\displaystyle \sum_{i-1}^{n_o{-1}} a_i\in A_\ell^{\prime\prime}\subseteq S^\prime, \displaystyle \sum_{i=1}^{n_o{-1}}a_1\leq\max S^\prime = 4J(n_o{-1},t)=$ order of $K.$ }
 we may choose a strictly increasing sequence $b_1,b_2,\ldots,b_{\theta_{n_o{-1}}}$ of elements of $K$, where 
$\theta_{n_o{-1}}=\displaystyle \sum_{i=1}^{n_o{-1}} a_i.$ In general define $\theta_i=\displaystyle \sum_{j=1}^i a_j,i=1,\ldots,n_o{-1}$ and define some subsets of $K$ as follows (see Fig. 1).\\
\begin{equation*}
\begin{split}
B_1&=\{b_1,\ldots,b_{a_1}\},\\
B_i&=\{b_{\theta_{i-1}+1},\ldots,b_{\theta_i}\}, ~~~~i=2,\ldots,n_o{-1},\\
P_1&=\left\{b_{\frac{a_1}{2}-1},\ldots,b_{a_1-1}\right\},\\
P_2&=\left\{b_{a_1+2},\ldots,b_{\frac{3a_1}{2}+2}\right\}, 
\end{split}
\end{equation*}
~~ and ~~~~~~~~~~~~~
$B_{n_o}=P_1\cup P_2.$

\[
\underbrace{\overbrace{b_1,\ldots,b_{\frac{a_1}{2}-1},\ldots,b_{a_1-1},b_{a_1}}^{\mathrm{B_1}}}_{\mathrm{P1}},~~~~
\underbrace{\overbrace{b_{a_1+1},b_{a_1+2},\ldots,b_{\frac{3a_1}{2}+2},\ldots,b_{a_1+a_2}}^{\mathrm{B_2}}}_{\mathrm{P_2}},~~~~\overbrace{b_{a_1+a_2+1},\ldots}^{B_3}
\]
~~~~~~~~~~~~~~~~~~~~~~~~~~~~~~~~~~~~~~~~~~Figure 1.
\\
 We claim that $f(B_i), i=1,\ldots,n_o$ and $A_\ell$ are the required numbers and class, i.e., 
$$~~~~~~~~~~~~\sum_{i\in I}f(B_i)\in A_\ell,~~~~ I \subseteq\{1,\ldots,n_o\}.~~~~~~~~~~~~~~~~~~~~~T$$
We first state two facts about the function $f$.
\begin{enumerate}
\item[$F_1$]~~$f(A\cup B)=f(A)+f(B)$ if $A$ and $B$ are of even order and $A<B\footnote{$A<B$ denotes here that the largest element of $A$ is less than the smallest element of $B$.}.$
\item[$F_2$]~~$f(A)=f(A-P)-f(P)$ when $P$ is a consecutive\footnote{ By a \underline{consecutive} \underline{subset} we mean that $a\in A,x<a<y,x,y\in P$ implies $a\in P.$} subset of $A$ of even order such that all the elements of $A$ larger than the largest element of $P$ are odd in number.
\end{enumerate}
\noindent
\par We first show $T$ assuming either $n_o\notin I$ or $n_o\in I$ but $1,2\notin I$. By $F_1$ we have $\displaystyle \sum_{i\in I}f(B_i)=f(\displaystyle \bigcup_{i\in I}B_i).$ But $\displaystyle \bigcup_{i\in I}B_i$ is of order $\displaystyle \sum_{i\in I}a_i$ if $n_o\notin I$ and of order $\displaystyle \sum_{i\in I-\{n_o\}}a_i+a_1$ if $n_o\in I,1,2\notin I$ since the $B_i$ are disjoint in either case and $|B_{n_o}|=a_1$ and $|B_i|=a_i,~~i=1,\ldots,n_o{-1}.$ Thus in either case $\displaystyle \bigcup_{i\in I}B_i\in A_{\ell}^\prime$
(since $|\displaystyle \bigcup_{i\in I}B_i|\in A_{\ell}^{\prime\prime}$) so $\displaystyle \sum_{i\in I}f(B_i)=f(\displaystyle \bigcup_{i\in I}B_i)\subseteq A_{\ell} $. Next we consider the case $n_o,1\in I,2\notin I$. We claim $f(B_{n_o})+f(B_1)=f((B_1-P_1)\bigcup P_2).$ [ This follows from the fact that $\frac{a_1}{2}$ is even so $f(B_{n_o})=f(P_1\bigcup P_2)=f(P_1)+f(P_2)$ (by $F_1$ since $|P_1|=|P_2|=\frac{a_1}{2}$ and $P_1<P_2$). 
Then by $F_2, f(B_1)+f(P_1)=f(B_1-P_1)$ and by $F_1, f(B_1-P_1)+f(P_2)=f((B_1-P_1)\bigcup P_2).$] But $(B_1-P_1)\bigcup P_2
$ is of order $a_1$ and $(B_1-P_1)\bigcup P_2< B_i, i=3,\ldots,n_o{-1}$,\\
so by $F_1,$~~~
 $  
\sum_{i\in I}f(B_i)=f(B_{n_0})+f(B_1)+\sum_{i\in I-\{n_0,1\}}f(B_i) =
$ 
$f((B_1-P_1)\cup\ P_2)+\sum_{i\in I-\{n_o,1\}}f(B_i)=f((B_1-P_1)\cup P_2\bigcup_{i\in I-\{n_o,1\}}B_i~.$
$\mbox{   But   } (B_1-P_1)\bigcup P_2\cup\bigcup_{i\in I-\{n_o,1\}}B_i \mbox{   is of order    }\sum_{i\in I-\{n_o\}}a_i.$
$ \mbox{ Thus } ((B_1-P_1)\bigcup P_2)\cup\bigcup_{i\in I-\{n_o,1\}}B_i \in A_{\ell},$$ ~\mbox{so}~\sum_{i\in I}f(B_i) =f(\bigcup_{i=1-\{1,n_o\}}B_i\cup (B_1-P_1)\cup P_2)\in A_{\ell}.$

Similar arguments hold in the cases $n_o,~	2\in I, ~1\notin I$ 
and $n_o,2, 1\in I$ noting $f(B_{n_o})+f(B_2)=f((B_2-P_2)\cup B_1)$ 
and $f(B_{n_o})+f(B_2)+f(B_1)=f((B_1\cup B_2)-B_{n_o}). \Box$
\section{\small FURTHER GENERALIZATIONS}
\noindent
\par It is natural to ask whether either Lemma 1 or Theorem 2 generalize in the following way:\\
\underline{Lemma $1^\prime$.}~~~~ Let all the non-empty finite subsets of a countably infinite set $S$ be divided into $t$ classes $A_1,\ldots,A_t$. Then there exists a countably infinite subset $K\subseteq S$ and an infinite sequence $i_1,i_2,\ldots,~~~~ i\leq i_j\leq t,~~ i=1,2,\ldots,$ such that all the $j-$subsets of $K$ are contained in $A_{i_j},~~ j=1,2,\ldots~~~~.$\\
\underline{Theorem $2^\prime$.}~~~~ Let the natural numbers be divided into $t$ classes $A_1,\ldots,A_t.$ Then there exists an infinite sequence $a_1,a_2,\ldots,$ and a number $\ell,  1\leq\ell\leq t,$ such that $\displaystyle \sum_{i\in I}a_i \in A_\ell,$ for all (non-empty) finite sets $I$ of natural numbers.
\par A counter example to Lemma $1^\prime$ is as follows. Divide all the finite (non-empty) subsets $S$ of the natural numbers $N$ into two classes $A_1$ and $A_2$ by the rule:
\begin{equation*}
\begin{split}
S\in A_1~~~~~~~~~~~~ \mbox{   iff   } ~~~~~~~~~~~~ |S| \in  S,\\
S\in A_2~~~~~~~~~~~~ \mbox{   iff   } ~~~~~~~~~~~~ |S| \notin  S.
\end{split}
\end{equation*}
Then $K$ contains $k-$subsets which contain $k$ and $k-$subsets which do not, i.e., subsets which belong to $A_1$ and ones which belong to $A_2$.
\par We note that if Theorem $2^\prime$ were true then using an argument similar to that of Corollary 1.1 we could assume the $a_1$ are distinct, however, at present the validity of Theorem $2^\prime$ is not known.
\\
\\
 \par  ~~~~~~~~~~~~~~~~~~~~~~~~~~~~~~~BIBLIOGRAPHY  \\
  1.~~~~1.	  I. Schur, Uber die Kongruenz $x^m$ + $y^m$ = $z^m$ (mod p), Jahresbericht der Deut. Mathematiker-Vereiningung, Band 25 (1916) pp. 114-117. 
 
\end{document}